\documentclass[a4paper,11pt]{article}

\usepackage{graphicx,psfrag}
\usepackage{amssymb,amsmath,amsthm}
\usepackage{enumerate}  
\usepackage{a4wide}
\usepackage{hyperref}

\usepackage{pstricks}
\usepackage{pst-pdf}
\usepackage{wrapfig}
\usepackage[font={small}]{caption}  

\usepackage{authblk}

\numberwithin{equation}{section}
\allowdisplaybreaks[4]

\newtheorem{proposition}{Proposition}[section]
\newtheorem{theorem}[proposition]{Theorem}

\theoremstyle{definition}

\newtheorem{remark}[proposition]{Remark}

\newcommand{\C}[1]{\mathbf{C}^{#1}}

\newcommand{\modulo}[1]{{\left|#1\right|}}

\newcommand{\reali}{{\mathbb{R}}}

\renewcommand{\epsilon}{\varepsilon}
\renewcommand{\phi}{\varphi}
\renewcommand{\theta}{\vartheta}

\renewcommand{\d}[1]{\mathinner{\mathrm{d}{#1}}}
\renewcommand{\div}{\mathop{\rm div}}
\DeclareMathOperator*\argmin{argmin}
\DeclareMathOperator*\argmax{argmax}

\newcommand{\tr}{\mathop{\rm tr}}

\begin{document}

\title{Numerical Algorithm for Optimal Control of Continuity Equations}

\author[1, 2]{\small Nikolay Pogodaev}

\affil[1]{\footnotesize Krasovskii Institute of Mathematics and Mechanics \authorcr 
Kovalevskay str., 16, Yekaterinburg, 620990, Russia}
\affil[2]{\footnotesize Matrosov Institute for System Dynamics and Control
  Theory\authorcr Lermontov str., 134, Irkutsk, 664033, Russia}

\date{}

\maketitle

\begin{abstract}

  \noindent An optimal control problem for the continuity equation is considered. The aim
  of a controller is to maximize the total mass within a target set at a given
  type moment. An iterative numerical algorithm for solving this problem is
  presented.

  \medskip

  \noindent\textit{2010~Mathematics Subject Classification: 49M05}

  \medskip

\end{abstract}

\section{Introduction}

Consider a \emph{mass} distributed on $\mathbb{R}^n$ that drifts along a controlled vector field $\mathbf{v}=\mathbf{v}(t,x,u)$.
The aim of the \emph{controller} is to bring as much mass as possible to a target set $A$ by a time 
moment $T$. 

Let us give  the precise  mathematical statement of the problem. Suppose that
$\rho=\rho(t,x)$ is the density of the distribution
and $u=u(t)$ is a strategy of the \emph{controller}.
Then, $\rho$ evolves in time according to the continuity equation
\begin{equation}
\label{eq:system}
\begin{cases}
\rho_t + \div{}\hspace{-2pt}_x \left(\mathbf{v}\left(t,x,u(t)\right)\rho\right)=0,\\
\rho(0,x) = \rho_0(x),
\end{cases}
\end{equation}
where $\rho_0$ denotes the initial density.
Our aim is to find a control $u$ that \textbf{maximizes} the following integral
\begin{equation}
\label{eq:cost}
J[u] = \int_{A}\rho(T,x)\d x.
\end{equation}
Typically, $u$ belongs to a set $\mathcal{U}$ of admissible controls.
Here we take the following one:
\begin{equation}
\label{eq:admissible}
\mathcal{U} = \left\{u(\cdot)\;\text{is measurable},\; u(t)\in U\; \text{ a.e. }\; t\in[0,T]\right\},
\end{equation}
where $U$ is a compact subset of $\mathbb{R}^m$.

In this paper we propose an iterative method for solving 
problem~\eqref{eq:system}--\eqref{eq:admissible},
which is based on the needle linearization algorithm for classical optimal
control problems~\cite{Srochko}. Given an initial guess $u^0$, the algorithm
produces a sequence of controls $u^k$ with the property $J[u^{k+1}]\geq J[u^k]$,
for all $k\in \mathbb{N}$.

A different approach for numerical solution
of~\eqref{eq:system}--\eqref{eq:admissible}
was proposed by S. Roy and A. Borz\`i in~\cite{Roy2017}. The authors used a specific discretization
of~\eqref{eq:system} to produce a finite dimensional optimization problem. It
seems difficult to compare the efficiency of both algorithms, because one
was tested for 2D and the other for 1D problems.


Finally, let us remark that problem~\eqref{eq:system}--\eqref{eq:admissible}
is equivalent to the following optimal control problem for an ensemble of
dynamical systems:
\begin{displaymath}
  \text{Maximize}\quad \int\rho_0(x)\d x\quad\text{subject to}\quad
  \begin{cases}
  \dot y = - \mathbf{v}(T-t,y,u(t)),\\
  y_0\in A.
  \end{cases}
\end{displaymath}
Indeed, instead of transporting the mass, one can transport the target $A$ in
reverse direction
aiming at the region that contains maximal mass.

\section{Preliminaries}

We begin this section by introducing basic notation and assumptions that will be used
throughout the paper. Next, we discuss a necessary optimality condition lying at the core
of the algorithm.

\subsection{Notation}

In what follows, $\Phi_{s,t}$ denotes the flow of a time-dependent vector field
$\mathbf{w}=\mathbf{w}(t,x)$, i.e., $\Phi_{s,t}(x) = y(t)$, where
$y(\cdot)$ is a solution to the Cauchy problem
\begin{displaymath}
  \begin{cases}
  \dot y(t) = \mathbf{w}\left(t,y(t)\right),\\
  y(s) = x.
  \end{cases}
\end{displaymath}
Given a set $A\subset\mathbb{R}^n$ and a time interval $[0,T]$,
we use the symbol $A^t$ for the image of $A$ under the map $\Phi_{T,t}$, i.e., $A^t
= \Phi_{T,t}(A)$. The Lebesgue measure on $\mathbb{R}$ is denoted by $\mathcal{L}^1$.

\subsection{Assumptions}

\begin{itemize}
\item The map $\mathbf v\colon [0,T]\times \reali^n\times U\to\reali^n$
is continuous. 

\item The map $x\mapsto \mathbf v (t,x,u)$ is twice continuously differentiable, 
for all $t\in [0,T]$ and $u\in U$. 

\item There exist positive
constants $L$, $C$ such that $\modulo{\mathbf v(t,x,u)-\mathbf v(t,x',u)}\leq L|x-x'|$
and 
$\modulo{\mathbf v(t,x,u)}\leq C\left(1+|x|\right)$,
for all $t\in [0,T]$, $u\in U$, and $x,x'\in\reali^n$.

\item The initial density $\rho_0$ is continuously differentiable.

\item The target set $A\subset \mathbb{R}^n$ is a compact tubular neighbourhood, i.e., $A$ is a compact set that
  can be expressed as a union of closed $n$-dimensional balls of a certain positive
  radius $r$.
\end{itemize}

In addition, to guarantee the existence of an optimal control
(see~\cite{Pogodaev2016} for details), we must assume that
\begin{itemize}
\item the vector field $\mathbf v$ takes the form
\[
  \mathbf v(t,x,u) = \mathbf{v}_0(t,x) + \sum_{i=1}^l \phi_i(t,u)\mathbf{v}_i(t,x),
\]
for some real-valued functions $\phi_i$, and the set
\[
\Phi(t,U) = 
\begin{pmatrix}
\phi_1(t,U)\\
\cdots\\
\phi_l(t,U)
\end{pmatrix}
\subset
\reali^l
\]
is convex.
\end{itemize}


\subsection{Necessary Optimality Condition}

The necessary optimality condition for problem~\eqref{eq:system}--\eqref{eq:admissible}
looks as follows:

\begin{theorem}[\cite{Pogodaev2016}]
  \label{thm:pmp}
Let $u$ be an optimal control for~\eqref{eq:system}--\eqref{eq:admissible} and $\rho$ be the corresponding trajectory with $\rho_0\in \C1(\mathbb{R}^n)$.
Then, for a.e. $t\in [0,T]$, we have 
\begin{displaymath}
\int_{\partial A^{t}} \rho(t,x)\, \mathbf v\left(t,x,u(t)\right)\cdot \mathbf{n}_{A^{t}}(x)\d {\sigma(x)}
= 
\min_{{\omega\in U}}\int_{\partial A^{t}} \rho(t,x)\, \mathbf v(t,x,{\omega})\cdot \mathbf{n}_{A^{t}}(x)\d {\sigma(x)}.
\end{displaymath}
Here $A^{t} = \Phi_{t,T}(A)$, where $\Phi$ is the {\em phase flow} of the
vector field $(t,x)\mapsto \mathbf{v}\left(t,x,u(t)\right)$,
$\mathbf{n}_{A^t}(x)$ is the measure theoretic {\em outer unit normal} to $A^t$ at $x$,
$\sigma$ is the $(n-1)$-dimensional {\em Hausdorff measure}.
\end{theorem}

Let $I\subseteq [0,T]$ be a measurable set of Lebesgue measure $\epsilon$. Given
two controls $u$ and $w$, we consider their mixture
\begin{equation}
  \label{eq:pcontrol}
u_{w,I}(t) = 
\begin{cases}
w(t), & t\in I,\\
u(t), & \text{otherwise}.
\end{cases}
\end{equation}

The proof of Theorem~\ref{thm:pmp} gives, as a byproduct, the following
increment formula
\begin{equation}
  \label{eq:increment}
J[u_{w,I}] - J[u] 
= \int_I 
\int_{\partial A^{t}} \rho(t,x)\, 
\left[
\mathbf v\left(t,x,u(t)\right)
- 
\mathbf v\left(t,x,w(t)\right)
\right]
\cdot \mathbf{n}_{A^{t}}(x)\d {\sigma(x)}
\d t
 + o(\epsilon),
\end{equation}
which will be used in the next section.

\section{Numerical Algorithm}

In this section we describe the algorithm, prove the improvement property
$J[u^{k+1}]\geq J[u^k]$, and discuss a possible implementation.

\subsection{Description}

\begin{enumerate}
  \item Let $u^k$ be a current guess. For each $t$, compute the set $\partial A^t$ and
    $\rho(t,\cdot)$ on $\partial A^t$.
  \item For each $t$, find 
    \begin{equation}
      \label{eq:first_max}
      w(t) = \argmin\left\{  \int_{\partial A^{t}} \rho(t,x)\, 
        \mathbf v\left(t,x,\omega\right)
        \cdot \mathbf{n}_{A^{t}}(x)\d {\sigma(x)}\;\colon\; \omega\in U\right\}.
      \end{equation}
  \item Let
    \begin{displaymath}
      g(t) = \int_{\partial A^{t}} \rho(t,x)\, 
      \left[
      \mathbf v\left(t,x,u^k(t)\right)
      - 
      \mathbf v\left(t,x,w(t)\right)
      \right]
      \cdot \mathbf{n}_{A^{t}}(x)\d {\sigma(x)}.
    \end{displaymath}
  
  \item For each $\epsilon\in (0,T]$, find
    \begin{equation}
      \label{eq:second_max}
      I(\epsilon) = \argmax\left\{  \int_{\iota} g(t)\d t\;\colon\;
       \iota\subset [0,T]\; \text{is measurable and}\;
        \mathcal{L}^1(\iota) = \epsilon
      \right\}.
    \end{equation}
  \item Construct $u_{w,I(\epsilon)}$ by~\eqref{eq:pcontrol}.
  \item Compute
    \begin{equation}
      \label{eq:max_last}
      \epsilon^* = \argmax\left\{ J[u_{w,I(\epsilon)}]\;\colon\; \epsilon\in (0,T] \right\}.
    \end{equation}

  \item Let $u^{k+1} = u_{w,I(\epsilon^*)}$.
\end{enumerate}

The algorithm produces an infinite sequence of admissible controls. Of
course, any its implementation should contain obvious modifications
that would cause the algorithm to stop after a finite number of iterations. Note that it may happen
that problems~\eqref{eq:second_max} and~\eqref{eq:max_last} admit no solution. In
this case $I(\epsilon)$ and $\epsilon^*$ must be taken so that the values of the
corresponding cost functions lie near the supremums.

\subsection{Justification}

If $u^k$ satisfies the optimality condition then we obviously get that
$u^{k+j}=u^k$, for all $j\in \mathbb{N}$. In particular, this means that
$J[u^{k+1}]=J[u^k]$.

If $u^k$ does not satisfy the optimality condition then $\int_{I(\epsilon)} g(t)\d t > 0$, 
for all small $\epsilon>0$. By the increment formula~\eqref{eq:increment}, we
have
\begin{displaymath}
  J[u_{w,I(\epsilon)}] - J[u^k] = \int_{I(\epsilon)} g(t)\d t + o(\epsilon).
\end{displaymath}
Since the integral from the right-hand side is positive for all small $\epsilon$, we
conclude that $J[u^{k+1}] = J[u_{w,I(\epsilon^*)}] > J[u^k]$, as desired.

\subsection{Implementation Details}

The method was implemented for 2D problems. 
All ODEs are solved by the Euler method. The
set $\partial A$ is approximated by a finite number of points.
Below we discuss in details all non-trivial steps of the algorithm.

\subsubsection*{Step 1}

In this step we must compute $\rho(t,x)$ for all $t$ and $x$
satisfying $x\in \partial A^{t}$. Recall that
\[
\rho(t,x) = \frac{\rho_0(y)}{\det D\Phi_{0,t}(y)},
\quad\mbox{where}\quad
y = \Phi_{t,0}(x).
\]
Using Jacobi's formula, we may write 
\[
\frac{d}{dt} \left(\det D\Phi_{0,t}(y)\right) = \left(\det D\Phi_{0,t}(y)\right)\cdot \tr\left[D\Phi_{0,t}(y)^{-1}\frac{d}{dt}D\Phi_{0,t}(y)\right].
\]
Meanwhile, by the definition of $\Phi$, we have  
\[
\frac{d}{dt}D\Phi_{0,t}(y) = D_x \mathbf{v}\left(t,\Phi_{0,t}(y),u(t)\right)\cdot D\Phi_{0,t}(y).
\]
Combining the above identities gives
\begin{equation*}
  \label{eq:ode_lin}
\frac{d}{dt} \left(\det D\Phi_{0,t}(y)\right) = \left(\det D\Phi_{0,t}(y)\right)\, \div \mathbf{v}\left(t,\Phi_{0,t}(y),u(t)\right).
\end{equation*}

Thus, computing of $\rho(t,x)$ requires solving two Cauchy problems, one for
finding $\Phi_{0,t}(y)$ and one for finding $\det D\Phi_{0,t}(y)$.

\subsubsection*{Step 2}

In general, the optimization problem~\eqref{eq:first_max} is nonlinear, which
makes it difficult. On the other hand, in many cases $U$ and
$\mathbf{v}$ enjoy the following
extra properties:
\begin{itemize}
\item the set $U$ is convex and
the vector field $\mathbf{v}$ is affine with respect to the control:
\[
\mathbf{v}(t,x,u) = \mathbf{v}_0(t,x) + \sum_{i=1}^m \mathbf{v}_i(t,x)\, u_i.
\]
\end{itemize}

Now~\eqref{eq:first_max} becomes a convex optimization problem, and thus
it can
be effectively solved.

\subsection*{Step 4}

The problem~\eqref{eq:second_max} seems difficult at first glance.
But note that it is equivalent to the following one:
\begin{equation}
  \label{eq:lambdaopt}
  \text{Minimize}\quad
  l(\lambda):=\left| \mathcal{L}^1\left(\left\{ t\;\colon\; g(t)\geq \lambda \right\} \right) - \epsilon\right|
  \quad\text{subject to}\quad \lambda\in [\min g,\max g].
\end{equation}
Indeed, if $\lambda_*$ solves~\eqref{eq:lambdaopt}, then the set
$I = \left\{ t\;\colon\; g(t)\geq \lambda_* \right\}$ solves the original
problem~\eqref{eq:second_max}. To find $\lambda_*$ numerically, we may take a finite
mesh on the interval $[\min g, \max g]$ and look for a node that gives the minimal
value to $l(\cdot)$.

\subsubsection*{Step 7}

In this step the cost 
\begin{displaymath}
 \int_A \rho(T,x)\d x = \int_{A^0}\rho_0(x)\d x
\end{displaymath}
must be computed.
To that end, we must know the whole set $A^0$, while on
the other steps of the algorithm we deal only with the boundaries of $A^t$.
It is interesting to note that, under the additional assumption that \begin{itemize}
\item the target set $A\subset \mathbb{R}^n$ is contractible and its boundary
  $\partial A$ is an $(n-1)$-dimensional smooth surface, 
the knowledge of $\partial A^0$ is enough for computing the cost.
\end{itemize}

Indeed, since the target $A=A^T$ is contractible, the set $A^0$ is contractible as well. 
Any differential form on a contractible set is exact~\cite{TuBook}. Hence
$\rho_0\d x = \d a$, for some $(n-1)$-dimensional differential form
$\alpha$. Now the Stokes theorem gives:
\begin{displaymath}
  \int_{A^0}\rho_0\d x = \int_{\partial A^0} \alpha.
\end{displaymath}

Let us compute $\alpha$ in the 2D case to illustrate this approach. 
We must find a form $\alpha = a_1\d x_1
+ a_2\d x_2$ such that $\d \alpha = \rho_0\d x$. The latter equation holds when
\begin{displaymath}
  \rho_0 = \frac{\partial a_2}{\partial x_1} - \frac{\partial a_1}{\partial x_2}.
\end{displaymath}
Hence, to get the desired $\alpha$, we may take
\begin{displaymath}
  a_1(x_1,x_2) = \int_0^{x_1}\rho_0(\xi,x_2)\d \xi + \rho_0(0,x_2),\quad a_2\equiv 0.
\end{displaymath}

\section{Examples}

This section describes several toy problems, which we used for testing the algorithm.

\subsection{Boat}

Consider a boat floating in the middle of a river at night. Since it is dark, the boatmen
cannot see any landmarks, and therefore are unsure about the boat's position.
They want to reach a river island at a certain time with highest probability. How should they act?

\begin{figure}[ht!]
\begin{center}
  \includegraphics[height=5cm]{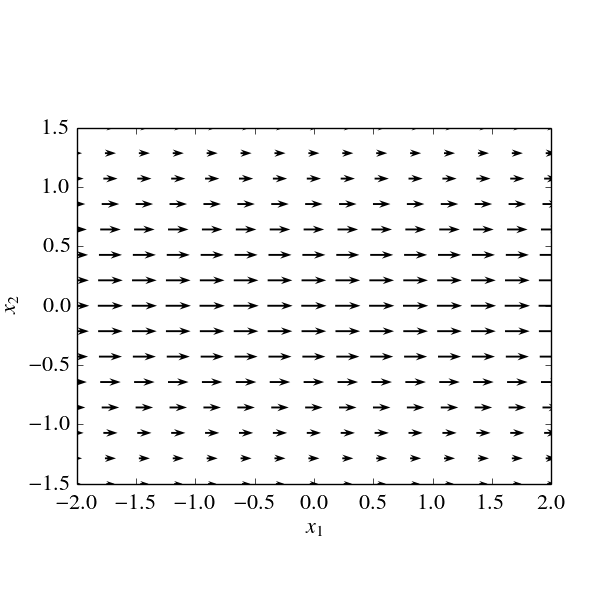}
  \qquad\qquad
  \includegraphics[height=5cm]{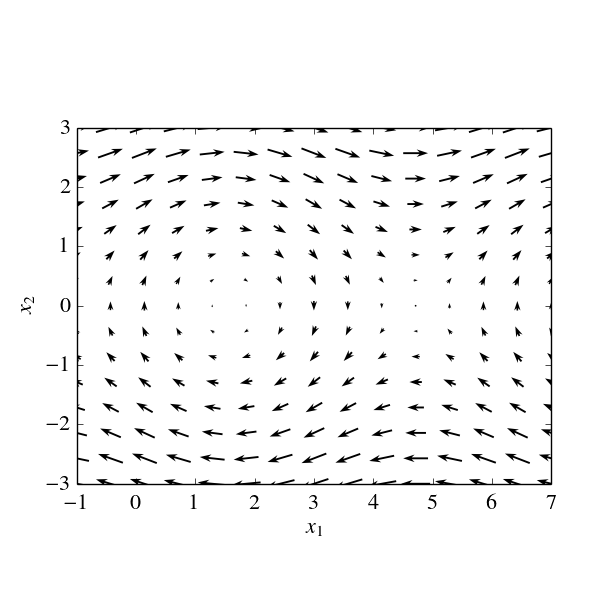}
\caption{Left: river drift. Right: pendulum drift. }
\label{fig1}
\end{center}
\end{figure}

\begin{figure}[ht!]
\begin{center}
  \includegraphics[height=3.5cm]{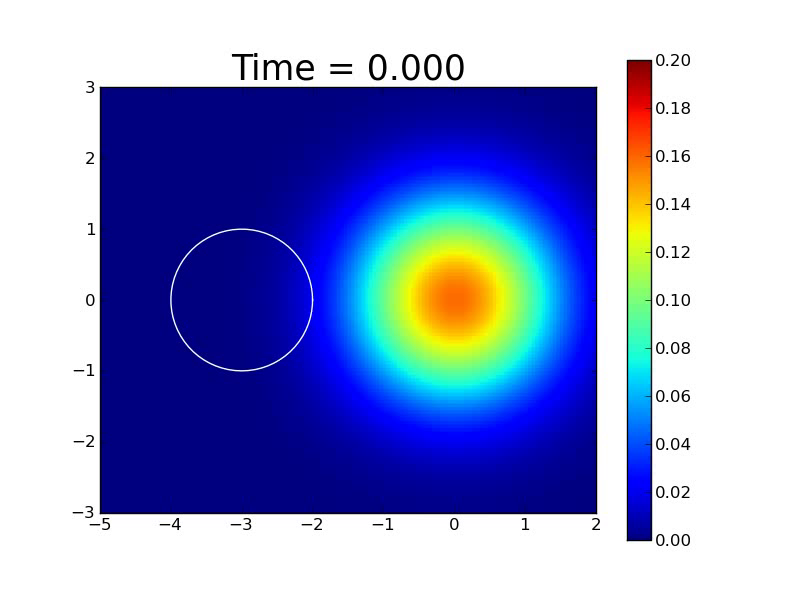} \quad
  \includegraphics[height=3.5cm]{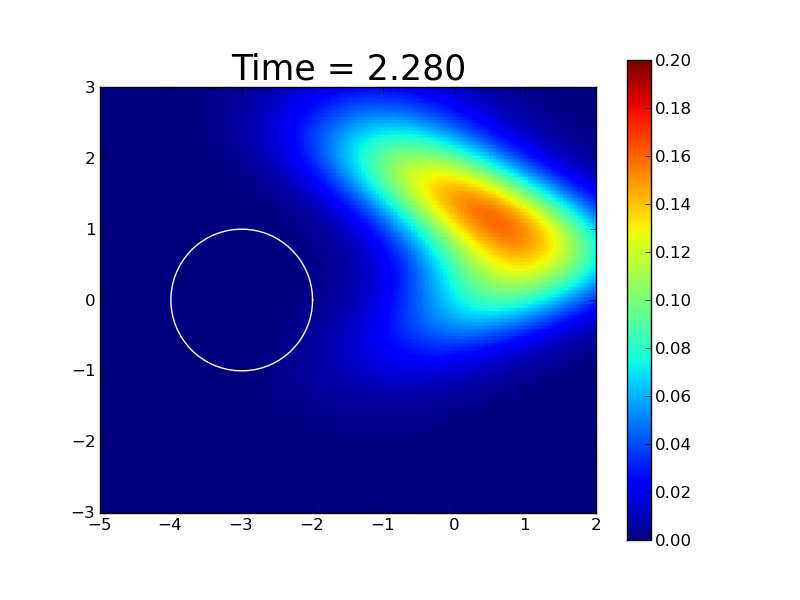} \quad
  \includegraphics[height=3.5cm]{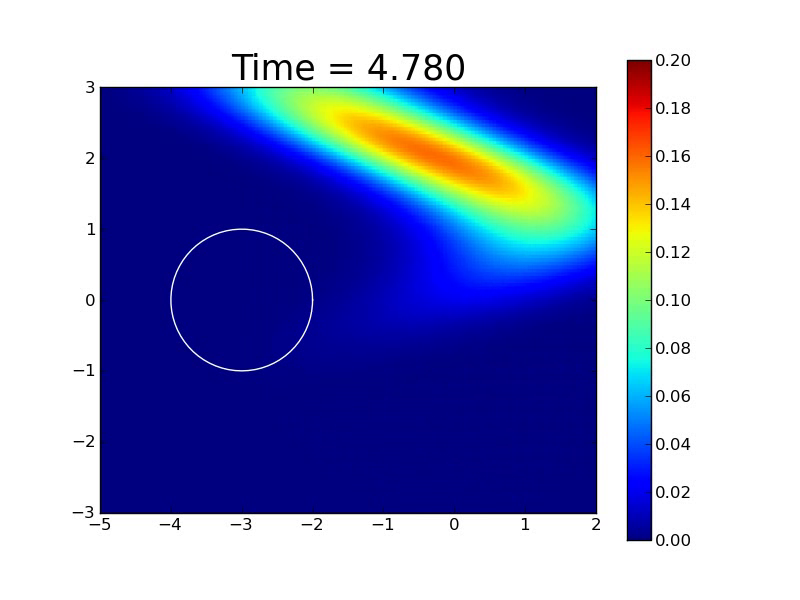} \quad
  \includegraphics[height=3.5cm]{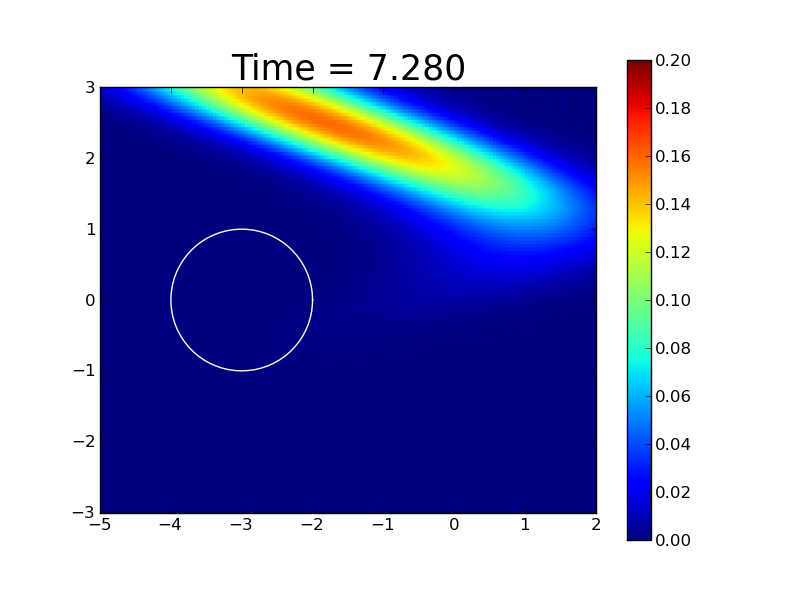} \quad
  \includegraphics[height=3.5cm]{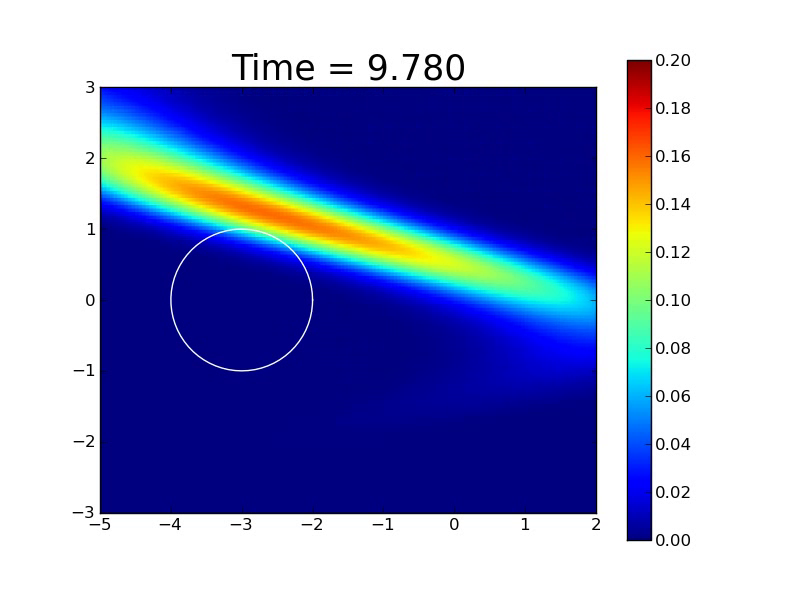} \quad
  \includegraphics[height=3.5cm]{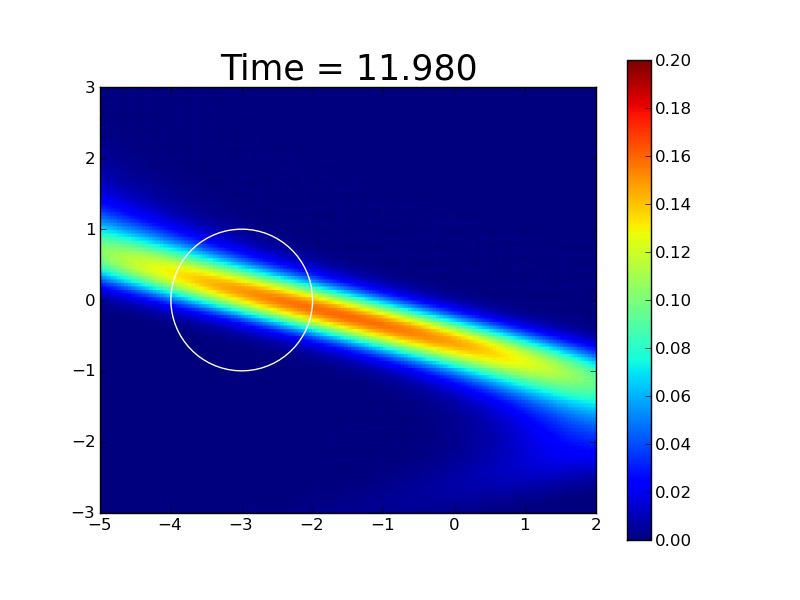} 
\caption{Trajectory for the boat problem computed by the algorithm.}
\end{center}
\end{figure}

Assume that the speed of the river water is given by
\[
\mathbf{v}_0(x) = 
\begin{pmatrix}
\alpha + e^{-\beta x_2^2}\\
0
\end{pmatrix},
\]
the island is a unit circle centered at $x_0$, the initial position of the boat is
described by the density function
\begin{equation}
\label{eq:normal}
\rho_0(x) = \frac{1}{2\pi\sigma^2}\, e^{-|x|^2/(2 \sigma^2)}.
\end{equation}
Thus, the boat's position $x(t)$ evolves according to the differential equation
\[
\dot x = \mathbf{v}_0(x) + u,
\]
where $u\in\mathbb{R}^2$ is a component of the boat's velocity due to rowing. Here $|u|\leq u_{\text{max}}$.

Parameters for the computation: $\sigma=1$, $\alpha=\beta=0.5$, $u_{\text{max}}=0.75$, $x_0=(-3,0)$, $T = 12$.

\subsection{Pendulum}

Here we want to stop a moving pendulum whose initial position is uncertain. In this case we have
\[
\mathbf{v}_0(x) = 
\begin{pmatrix}
x_2\\
\cos x_1
\end{pmatrix},
\qquad 
\mathbf{v}_1(x) = 
\begin{pmatrix}
1\\
0
\end{pmatrix}.
\]
Hence the control system takes the form
\[
\dot x = \mathbf{v}_0(x) + u\, \mathbf{v}_1(x),
\]
where $u\in [-u_{\text{max}},u_{\text{max}}]$ is an external force. The initial position of the pendulum is given by~\eqref{eq:normal}.
The target is a unit circle centered at $(\pi/2,0)$.

Parameters for the computation: $\sigma=1$, $u_{\text{max}}=0.5$, $x_0=(\pi/2,0)$, $T = 6$.
\begin{figure}[ht!]
\begin{center}
  \includegraphics[height=3.5cm]{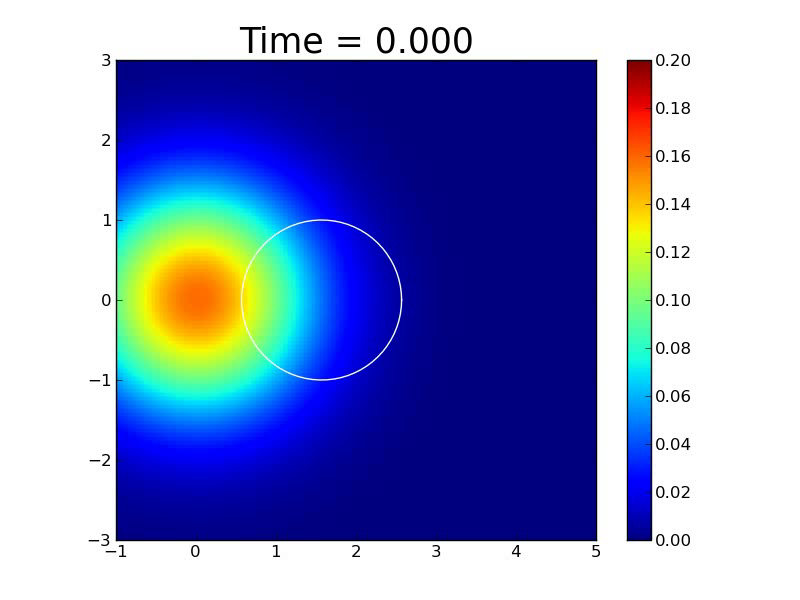} \quad
  \includegraphics[height=3.5cm]{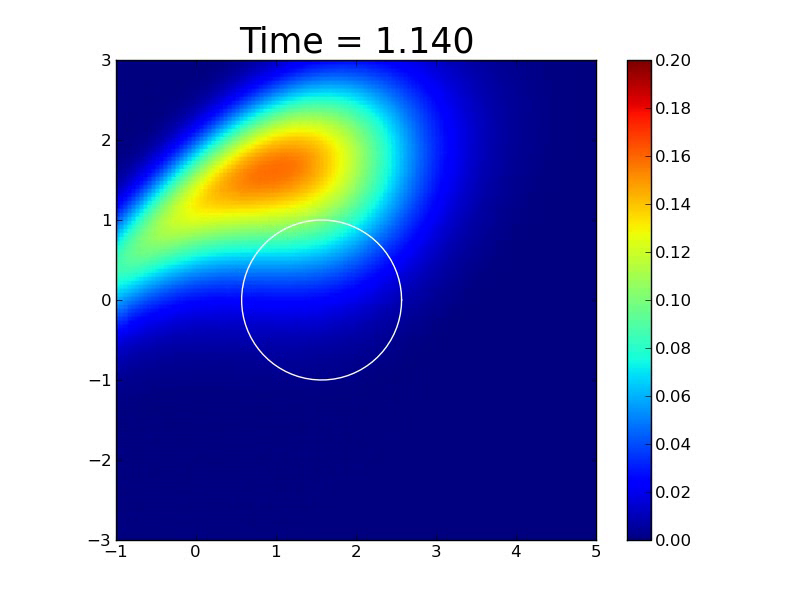} \quad
  \includegraphics[height=3.5cm]{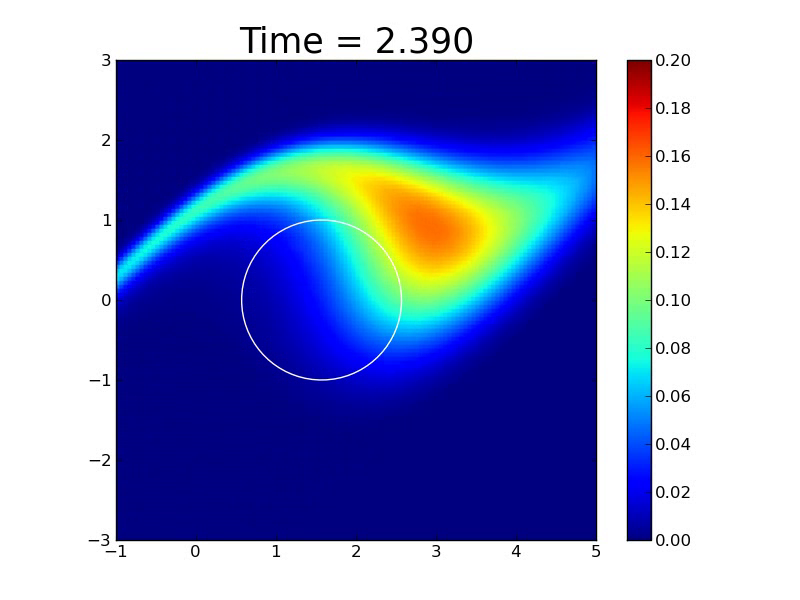} \quad
  \includegraphics[height=3.5cm]{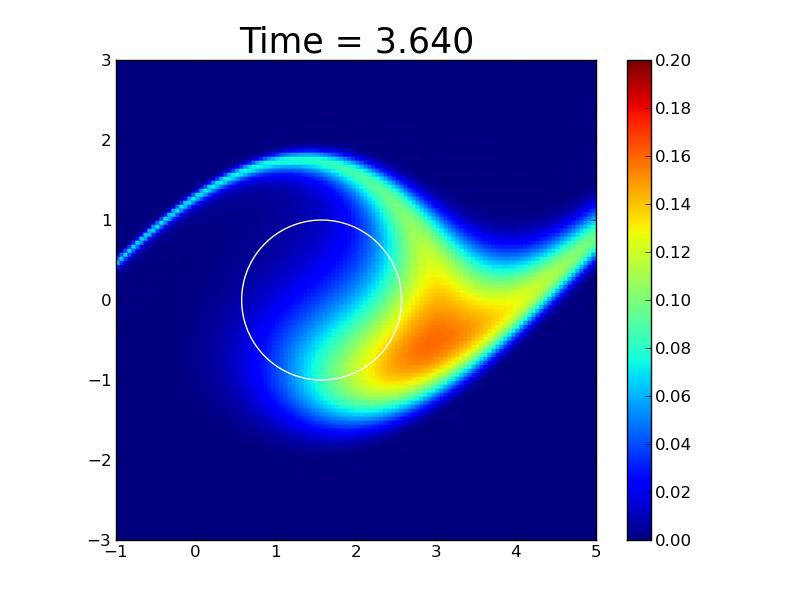} \quad
  \includegraphics[height=3.5cm]{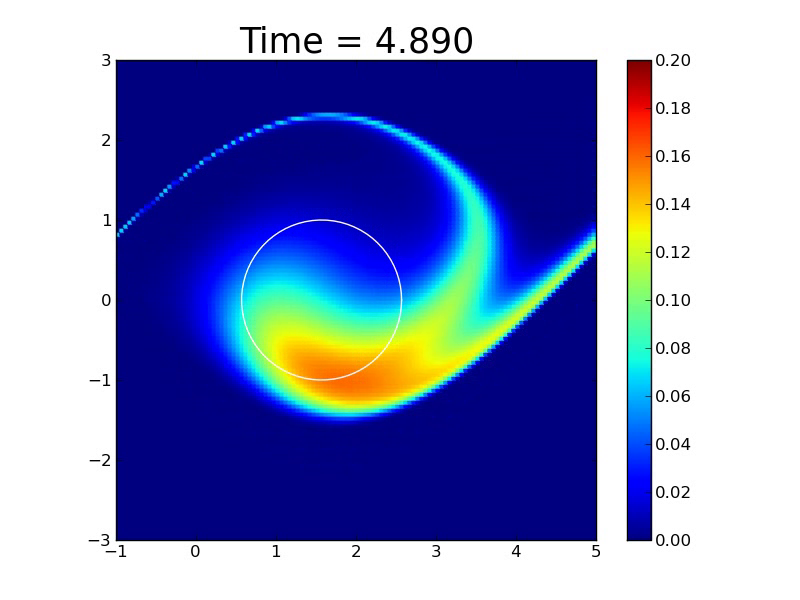} \quad
  \includegraphics[height=3.5cm]{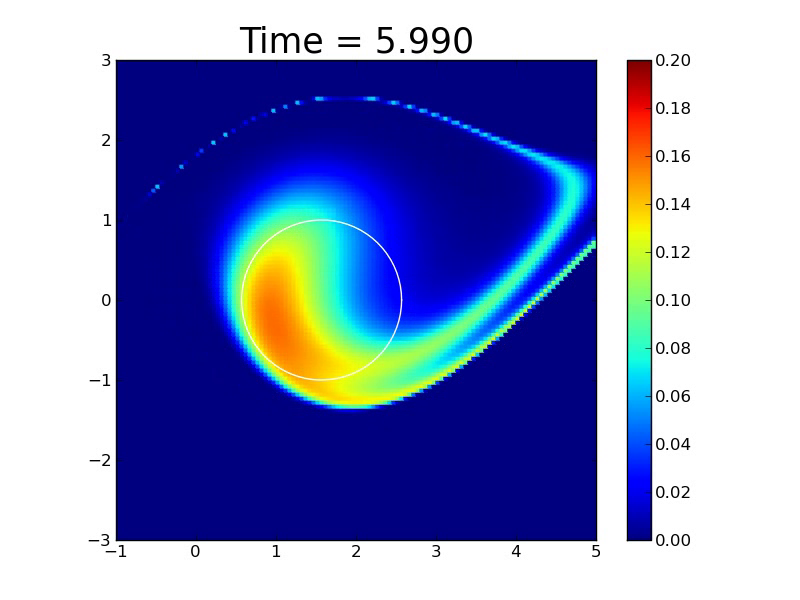} 
\caption{Trajectory for the pendulum problem computed by the algorithm.}
\end{center}
\end{figure}
\subsection{Sheep}

Consider a herd of sheep located near the origin. The sheep are effected by a vector field 
$v_0(x)$ pushing them away from the origin. 
To prevent this we can turn on repellers, which are located at the following positions
\[
x_k = \left(R\cos\frac{2\pi (k-1)}{m}, R\sin\frac{2\pi (k-1)}{m}\right),\qquad k=1,\ldots,m.
\]

\begin{figure}[ht!]
\begin{center}
  \includegraphics[height=6cm]{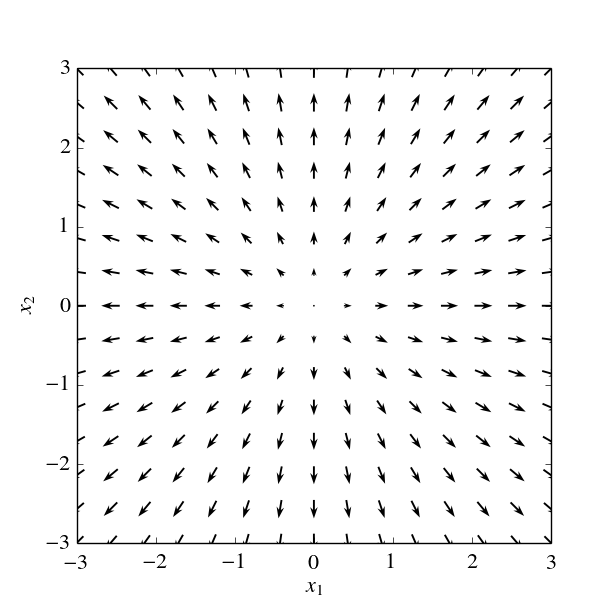}
  \qquad\qquad
  \includegraphics[height=6cm]{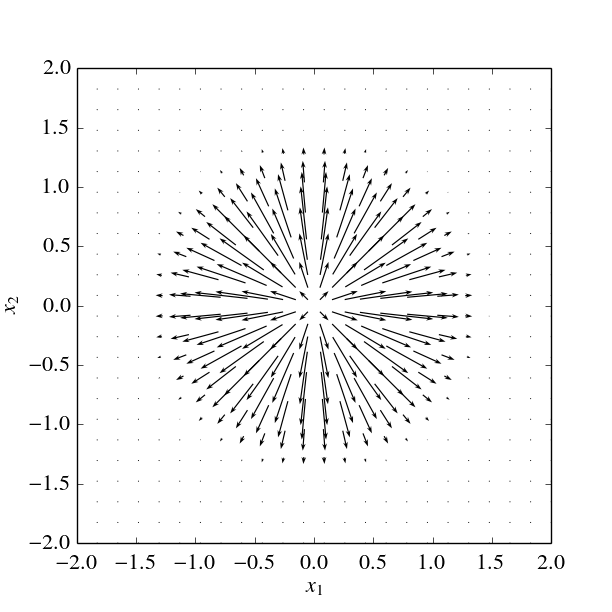}
\caption{Left: sheep drift. Right: repeller's force field. }
\label{fig2}
\end{center}
\end{figure}

\begin{figure}[ht!]
\begin{center}
  \includegraphics[height=3.5cm]{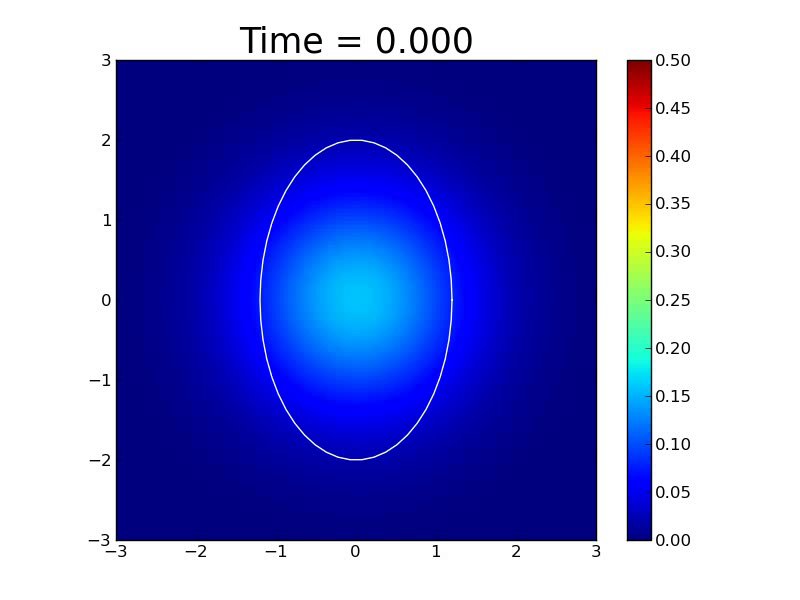} \quad
  \includegraphics[height=3.5cm]{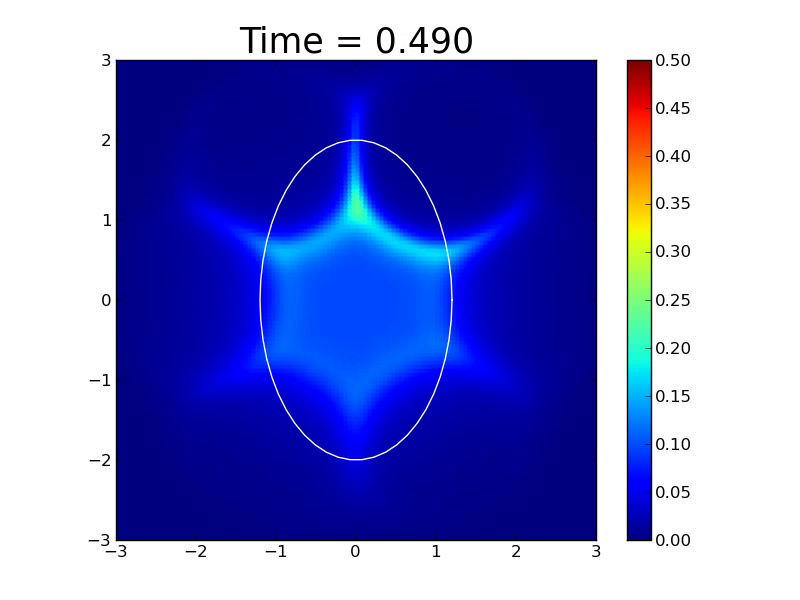} \quad
  \includegraphics[height=3.5cm]{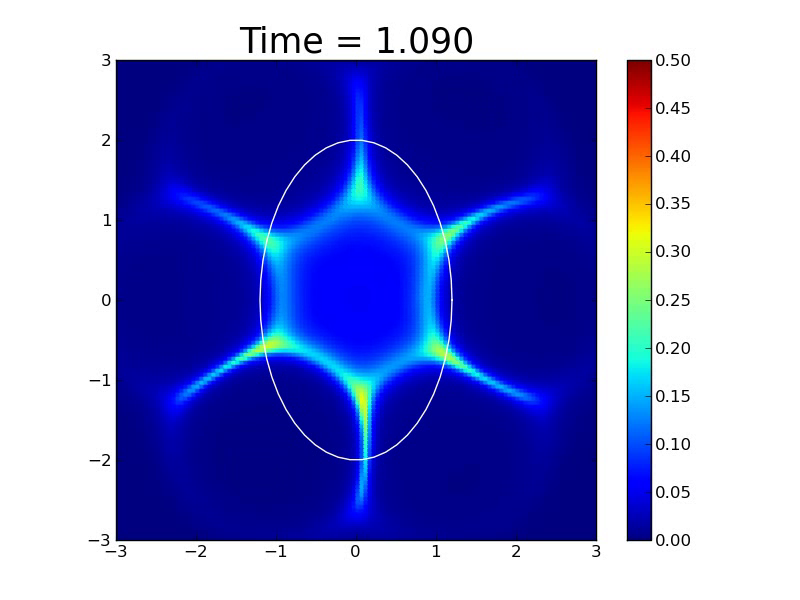} \quad
  \includegraphics[height=3.5cm]{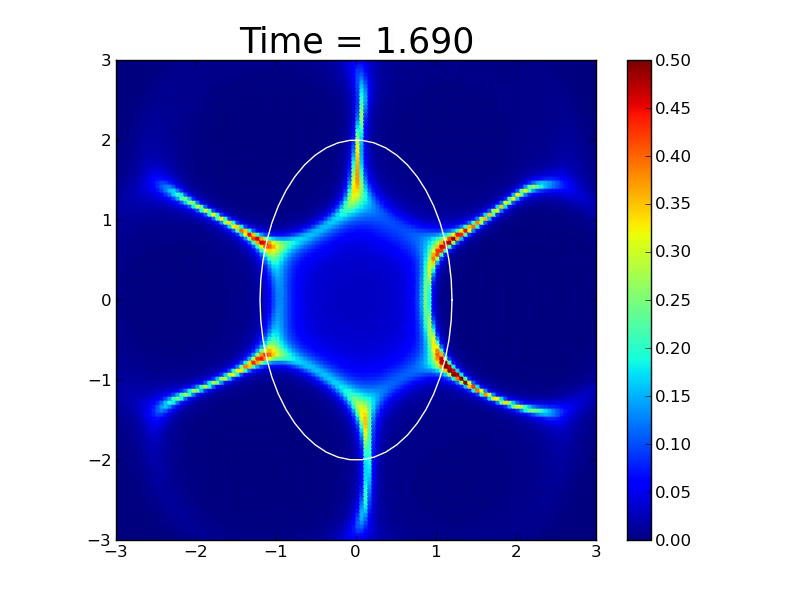} \quad
  \includegraphics[height=3.5cm]{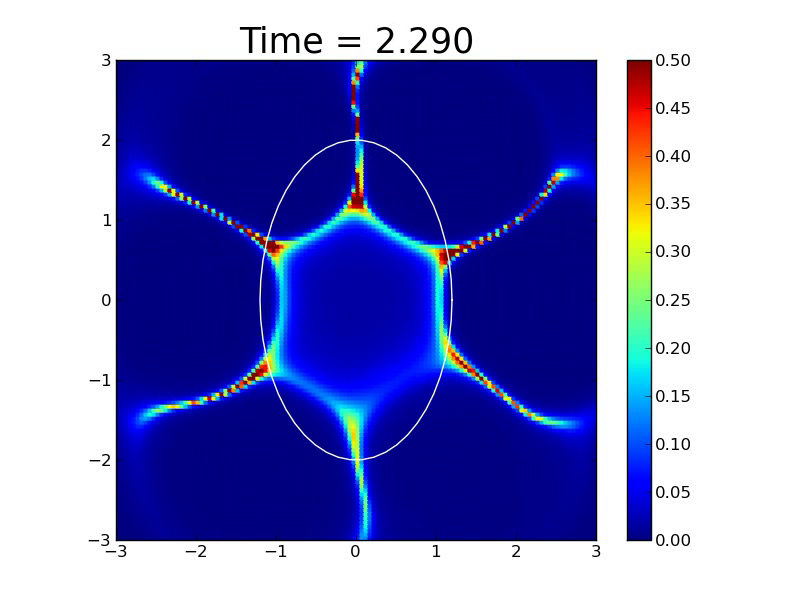} \quad
  \includegraphics[height=3.5cm]{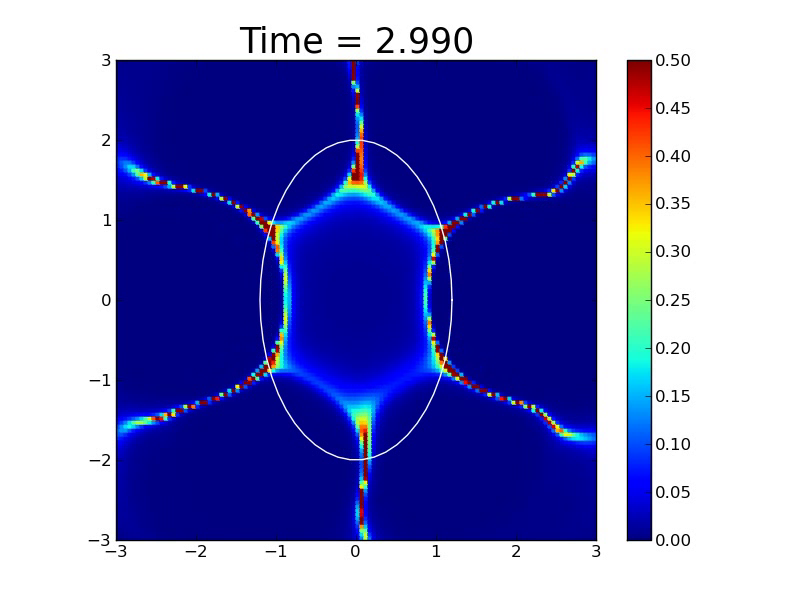} 
\caption{Trajectory for the sheep problem computed by the algorithm.}
\end{center}
\end{figure}

Each repeller produces a vector field $\mathbf{v}_k(x)$. So we have
\[
\mathbf{v}(x,u) = \mathbf{v}_0(x) + \sum_{k=1}^m u_k\mathbf{v}_k(x),
\]
where $u_k$ is an intensity of $k$-th repeller.
The control $u = (u_1,\ldots,u_m)$ belongs to the simplex
\[
U = \left\{(u_1,\ldots,u_m)\;\colon\;\sum_{k=1}^m u_k = 1,\, u_k\in [0,1],\; k=1,\ldots,m\right\}.
\]

In what follows we set
\[
\mathbf{v}_0(x) = \alpha\,\frac{x-x_0}{\sqrt{1+|x-x_0|^2}},
\]
where $x_0$ is a certain point not far from the origin, and
\[
\mathbf{v}_k(x) = \beta\, e^{-|x-x_k|^4}(x-x_k),\qquad k=1,\ldots,m.
\]

Suppose that the initial distribution is given by~\eqref{eq:normal}, the target is an ellipse centered at $x_0$ whose
major and minor semi-axes are $a$ and $b$.

Parameters for the computation: $\sigma=1$, $x_0=(0,0)$, $T = 3$, $m=6$, $a =
2$, $b = 1.2$.

\begin{remark}
The answer to the minimization problem
\[
\sum_{i=1}^m c_i \omega_i \to \min,\quad \omega\in U, 
\]
arising in the second step of the algorithm, is very simple. Let $j$ be such that
\[
c_j \leq c_i\quad\mbox{for all }i = 1,\ldots,m;
\]
then an optimal solution is given by $\bar\omega = (0,\ldots,0,1,0,\ldots,0)$, where $1$ is located at the $j$-th position. In particular, this
means that at every time moment $t$ only one repeller is turned on. Hence instead of repellers, we may think of a dog
that jumps from one place to another.
\end{remark}


\subsubsection*{Acknowledgements}

The work was supported by the Russian Science Foundation, grant No 17-11-01093.

\small{

  \bibliography{references}

  \bibliographystyle{abbrv}

}

\end{document}